\documentclass[12pt]{article}

\usepackage{amsmath}
\usepackage{graphicx}
\usepackage{amsfonts}
\usepackage{amsopn}
\usepackage{amsthm}

\newcommand{\zt}{\tilde{z}}
\newcommand{\rt}{\hat{\rho}}

\newcommand{\dm}{d\sigma}

\newcommand{\Db}{{\bar{D}}}

\newcommand{\Adm}{{{\rm{Adm}}(\Db,\Tb)}}
\newcommand{\Admr}{{{\rm{Adm}_{refl}}(\Db,\Tb)}}

\renewcommand{\r}{\rho}

\newcommand{\Cp}{C_{>0}(\bar{D})}

\renewcommand{\O}{\mathcal{O}}
\newcommand{\Tb}{{\bar{T}}}

\newcommand{\ra}{\rangle}
\newcommand{\la}{\langle}

\newcommand{\comment}[1]{}

\newcommand{\e}{\epsilon}
\renewcommand{\d}{\delta}

\newcommand{\z}{\zeta}

\newcommand{\F}{\mathcal{F}}

\renewcommand{\epsilon}{\varepsilon}

\renewcommand{\(}{\left(}
\renewcommand{\)}{\right)}

\DeclareMathOperator{\supp}{supp}

\DeclareMathOperator{\diam}{diam}

\newtheorem{theorem}{\bf \normalsize \bf  Theorem}[section]

\newtheorem{corollary}[theorem]{\bf \normalsize \bf Corollary}
\newtheorem{lemma}[theorem]{\bf \normalsize \bf Lemma}
\newtheorem{proposition}[theorem]{\bf \normalsize \bf Proposition}

\newtheorem{problem}{\bf \normalsize \bf  Problem}
\newtheorem*{problem'}{\bf \normalsize \bf  Problem I$_{\text{\bf w}}$}

\theoremstyle{remark}
\newtheorem{remark}[theorem]{\bf \normalsize \bf Remark}

\theoremstyle{definition}
\newtheorem{definition}[theorem]{\bf \normalsize \bf Definition}

\date{}
\title{A rigorous analysis using optimal transport theory
for a two-reflector design 
problem with a point source 
}
\author{Tilmann Glimm\\
Department of Mathematics\\
Western Washington University, Bellingham, WA 98225\\
}
\date{August 17, 2009}

\begin{document}


\maketitle
\thanks{
reflector design, geometric optics, optimal transportation}

\begin{abstract}
We consider the following geometric optics problem:
Construct a system of two reflectors which transforms
a spherical wavefront generated by a point source 
into a beam of parallel rays. This beam has a prescribed
intensity distribution. We give a rigorous
analysis of this problem. The reflectors we construct are 
(parts of) the boundaries of convex sets. We prove existence
of solutions for a large class of input data and give
a uniqueness result. To the author's knowledge, this
is the first time that a rigorous mathematical analysis of this
problem is given. The approach is based on optimal transportation
theory. It yields a practical algorithm for finding the reflectors.
Namely, the problem is equivalent to a constrained linear optimization
problem.
\end{abstract}
\section{Introduction}
We consider the following beam shaping problem from geometric optics, 
sketched in Figure~\ref{fig1}:
Suppose we are given a spherical wavefront with
a given intensity distribution
emitted from a point source. We would like to transform this
input beam into an output beam of parallel light rays with a desired intensity distribution.
This transformation is to be achieved with a system of two reflectors\footnote{A treatment with three reflectors 
is possible, see \cite{oliker95}. 
For reasons of space limitations, a two reflector
system is sometimes more desirable from an engineering point of view.}.
This paper deals with the mathematical problem of finding 
these reflecting surfaces for given input and output apertures and
input and output intensities\footnote{From a practical
point of view, an additional requirement is that no blockage of the light 
may occur. We do not deal with this explicitly, but we note that this can 
generally be achieved by choosing the reduced
optical path length  large enough. (See Section~\ref{sect_statement}
for the definition of the reduced optical path length.)}.
The geometric optics approximation is assumed.

\begin{figure}
\begin{center}
\end{center}
\caption{Geometry of the reflector problem. Note our convention regarding spatial coordinates, illustrated 
by the coordinate system in the lower left hand corner: The output beam propagates in the direction
of the {\em{negative}} $z-$axis, and points in the plane perpendicular to the $z-$axis are denoted
by the vector $x\in\mathbb R^2$. Thus a generic point in three dimensional space is denoted by $(x,z)$. 
For more details see the accompanying text in 
Section~\ref{sect_statement}.}\label{fig1}
\end{figure}

These types of problems are of practical interest and there exists an 
extensive engineering literature on them;
see Section 5 in V. Oliker's stimulating survey paper \cite{OlikerSurvey} and the many sources cited there. As explained there,
the approaches described in the engineering literature
are usually only applicable to specific data and are usually only
justified by specific numerical examples. 
In \cite{OlikerSurvey}, the problem was reformulated as
equation of Monge-Amp\`ere type for the polar radius
of the first reflector. A rigorous mathematical analysis,
showing the existence of solutions,
was however lacking as stated in \cite{OlikerSurvey}.  

In the present paper, we provide this rigorous
mathematical analysis, using a novel approach
to the problem. (To the best of our knowledge, this is the
first time such a rigorous analysis has been provided.)
We show the existence of solutions for given input and 
output intensities and provide a uniqueness result for the
ray tracing map associated with the reflector system.
The two reflectors we obtain
are always concave; they are contained in the boundary
of certain convex sets. 
 The approach is based on 
advances in the theory of Monge-Amp\`ere equations and
optimal transportation that have been made in about the last 15 years
\cite{Brenier1,Caf_alloc:96, Gangbo/Mccann:95,villani}. 
Similar approaches have fruitfully previously been applied to other beam shaping
 problems by X-J Wang in \cite{wang}, and
independently by V. Oliker and the author in  \cite{go04,go03}. 

We now describe the approach and results in more detail. For this,
denote the input aperture by $\Db$, and the output aperture by $\Tb.$
Thus $\Db$ is a subset of the unit sphere $S^2$, and $\Tb$ is contained
in a plane perpendicular to the direction of the output beam. (See again
Figure~\ref{fig1}.)
The input and output intensities are given by
nonnegative functions $I(m), m\in\Db$,
and $L(x)$, $x\in\Tb$. 

We represent the two reflectors via the polar radius $\r(m), m\in\Db$
for the first one, and as the graph
of a function $z(x)$, $x\in\Tb$ for the second one. See again Figure~\ref{fig1}.  
One of the main results is that finding the functions 
 $\r(m)$ and $z(x)$ is equivalent to solving the following
constrained optimization problem for certain transforms
of $\r(m)$ and $z(x)$:
\begin{align}
&\text{Minimize }\int_\Db \log\rt(m)I(m) \dm+\int_\Tb \log\zt(x)L(x)dx\label{constr_opt1}\\
&\text{subject to }\log\rt(m)+\log\zt(x)\geq \log K(m,x){\text{ for all }}m\in\Db,x\in\Tb.\label{constr_opt2}
\end{align}
Here $\log\rt(m)$ and $\log\zt(x)$ are certain transforms of $\r(m)$ and $z(x)$, respectively,
given explicitly in Definitions~\ref{def_zt} and \ref{def_rt} below. The constraint is given   by
\[
 \log K(m,x)=\log\left[\frac{\ell - \la m_x,x\ra}{2\ell(\ell^2-|x|^2)(1+m_z)}-\frac{1}{4\ell^2}\right].
\]
Here $m=(m_x,m_z)\in\Db\subseteq S^2$, $x\in\Tb\subseteq\mathbb R^2$, 
and $\ell>0$ is a constant. For details see the text below.

Note that this reformulation not only gives a theoretical existence result for solutions, but
it also translates into a practical method for finding the solution. In fact,
the discretization of the constrained optimization problem~(\ref{constr_opt1})-(\ref{constr_opt2})
is a standard linear programming problem and can be solved numerically.

As mentioned before, the approach of this paper has a strong connection to, and is in fact motivated by,
the theory of optimal transportation. (See for example \cite{Brenier1,Caf_alloc:96, Gangbo/Mccann:95},
and in general the recent survey  \cite{villani} and the extremely 
extensive bibliography cited there.) Consider the ray tracing map, or reflector map, $\gamma\colon\Db\to\Tb$.
So a ray emitted in the direction $m\in\Db$ will be transformed by the reflector
system into a ray labeled by $x=\gamma(m)\in\Tb$. (See again Figure~\ref{fig1}, and also Figure~\ref{refl_sketch}.)

Consider the optimal transportation problem for moving the
measure $I(m)\dm$ on $\Db\subseteq S^2$ to the measure $L(x)dx$ on $\Tb\subseteq\mathbb R^2$
via a transformation $P\colon\Db\to\Tb$
in such a way that the ``transportation cost'' 
\[
    C(P)=\int_{\Db}\log K(m,P(m))I(m)\dm   
\]
is maximized\footnote{It is unclear to the author whether
the cost $C(P)$ has any physical meaning. If {\em{maximizing}} a transportation cost seems unintuitive,
an equivalent formulation is of course to {\em{minimize}} $-C(P)$}. We prove that the ray tracing map $\gamma$ actually 
solves this problem. This gives a uniqueness result for the reflector map, 
see Theorem \ref{thm_ex_unique}.

The structure of the paper is as follows: In Section~\ref{sect_statement}, we state the
problem and our notation precisely and introduce the transforms $\rt(m)$
and $\zt(m)$. In the next section, we give certain analytical definitions of our concepts of reflectors
and the corresponding reflector map. In Section~\ref{section_geometry}, we then 
show that these analytical constructions correspond to some geometric constructions.
This geometric content is that we seek to represent the first reflector
as the boundary of the intersection of a certain family of spheroids, and we seek to represent the second
reflector as the boundary of the intersection of a family of paraboloids. 
We also justify that our abstract definition of the reflector map $\gamma$
is consistent with the optical definition.
In Section~\ref{section_equiv}, we formulate the
reflector problem and show that it is equivalent to the constrained optimization
problem (\ref{constr_opt1})-(\ref{constr_opt2}). Then in the next section, we prove the existence of solutions,
which follows from a standard compactness argument. Finally, in Section~\ref{transp}, we show the connection
to the transportation problem mentioned above, and finally state the main theorem, Theorem~\ref{thm_ex_unique}, 
on the existence and
uniqueness of solutions to the reflector problem.

Throughout, the treatment is very similar to the papers \cite{go04,go03}, which treated
similar but distinct beam shaping problems.  The main difference is that the problem at hand
requires the introduction of the transforms $\rt(m)$ and $\zt(x)$, and a more complicated form
for the cost function $\log K(m,x)$. These complications arise mostly because the first reflector
is best described using polar coordinates, and the second reflector is best described using
Cartesian coordinates. For the sake of completeness and being self-contained, 
we include most of the proofs here, leaving out some details if they can easily be filled in from
\cite{go04} and \cite{go03}.

\section{Statement of the problem and assumptions}\label{sect_statement}
We first fix our notations and assumptions in this section.
Consider the configuration show in Figure~\ref{fig1}. A point source located at the origin
$O=(0,0,0)$ generates a spherical wave front over a given {\em{input aperture}} $\Db$ contained in the
unit sphere $S^2$. The latter is required to be transformed into a beam of parallel rays
propagating in the direction of the negative $z-$axis by means of a system of two reflectors.
A cross section of the output beam is specified as a set on a plane perpendicular to the direction of
propagation. Such a cross section is called the {\em{output aperture}}, and denoted by $\Tb$. 

We denote points in space $\mathbb R^3$ by pairs $(x,z)$, where $x\in \mathbb R^2$ is the position vector 
in a plane perpendicular to the direction of propagation and $z\in\mathbb R$ is the coordinate in the
(negative) direction of propagation. See again Figure~\ref{fig1} for our convention on the
direction of the $z-$axis. Points on the unit sphere $S^2$ will typically be denoted by $m\in S^2$;
their components are also written as $m=(m_x,m_z)$ with $|m_x|^2+m_z^2=1$.

We fix the output aperture in the plane $z=-d$. We will seek to represent the two reflectors
as the graph of the polar radius $\r(m)$ and as the graph of a function $z(x)$, respectively, 
as shown in Figure~\ref{fig1}. That is
\begin{align*}
& \text{Reflector 1: }\Gamma_\r=\{\r(m)\cdot m\;\bigl|\; m\in \Db\},\\
& \text{Reflector 2: }\Gamma_z=\{(x,z(x))\;\bigl|\; x\in \Tb\}.
\end{align*}

We now have the following assumptions on the input data. Let $D$ be an open, nonempty subset of
$S^2$ with closure $\Db$ such that $(0,0,1)\notin \Db$ and $(0,0,-1)\notin \Db$. So  there is some $\e>0$ 
such that for $m=(m_x,m_z)\in\Db$, 
we have $1-\e>m_z>-1+\e$. Let further $T$ to be an open, bounded, nonempty subset of
$\mathbb R^2$. Denote its closure by $\Tb$.

The geometrical optics approximation is assumed. It follows from general principles
of geometric optics that  all rays will have equal length from $(0,0,0)$
to the plane $z=-d$; this length is called the optical path length and will be denoted by $L$.
We define the {\em{reduced optical path length}} as $\ell=L-d$. In terms of the mathematical
problem, $\ell$ is an input parameter.
We pick the reduced optical path length $\ell>0$ large enough so that the following
conditions are satisfied:
\begin{align}
&\ell> \max_{x\in\Tb}|x|,\quad\quad {\text{ and }}\\
& \frac{\ell - \la m_x,x\ra}{(\ell^2-|x|^2)(1+m_z)}-\frac{1}{2\ell}>0 {\text{ for all }}m\in\Db, x\in\Tb.
\label{def_K}
\end{align}
(It is not hard to see that the second condition is indeed satisfied for large enough $\ell$.) 

We'll also use the convenient notation 
\[
\delta =\frac{1}{2\ell}.
\]
Finally, we define the following two transformations, which are central
to the analysis:
\begin{definition}\label{def_zt}
Let $z=z(x)$ be a continuous function defined on $\Tb\subseteq\mathbb R^2$.
Then define the function
\begin{equation}
\zt(x)=\delta - \frac{z(x)}{\ell^2-|x|^2}\quad\quad \text{   for }x\in\Tb.
\end{equation}
\end{definition}

\begin{definition}\label{def_rt}
Let $\r=\r(m)$ be a continuous function defined on $\Db\subseteq S^2$ with $\r>0$.
Then define the function
\begin{equation}
\rt(m)=-\delta + \frac{1}{2\r(m)\cdot (m_z+1)}\quad\quad \text{   for }m\in\Db.
\end{equation}
\end{definition}
These transformations are obviously invertible, namely
\begin{align*}
& z(x)=(-\zt(x)+\d)(\ell^2-|x|^2)\quad\quad  \text{   for }x\in\Tb,\\
& \r(m)=(2(m_z+1)(\rt(m)+\d))^{-1} \quad\quad \text{   for }m\in\Db. 
\end{align*}

The following lemma is now obvious with the above formulas.
\begin{lemma}
\begin{enumerate}
\item The transformation 
$
C(\Tb)\to C(\Tb), z\mapsto \zt, 
$
is a bijection.

\item Denote by $\Cp$ the set of all positive continuous functions on $\Db$. Then the transformation
$
\Cp \to  \{\sigma\in C(\Db)\;\bigl|\; \sigma>-\delta\}, \r \mapsto \rt,
$
is a bijection.
\end{enumerate}
\end{lemma}

\section{Reflector pairs and the reflector map}
We now come to the central definition of this paper, namely that of a reflector pair.
We first give the analytic definition. In the next section, we look at the 
geometric interpretation. In preparation, we first define the following function:
\begin{definition}
Define the function
\[
   K(m,x)=\delta\frac{\ell - \la m_x,x\ra}{(\ell^2-|x|^2)(1+m_z)}-\delta^2 {\text{ for }}m=(m_x,m_z)\in\Db, x\in\Tb.
\]
Note that $K(m,x)>0$ by (\ref{def_K}). The function $\log K(m,x)$
can be interpreted as some kind of
cost of transporting a unit of energy from $m$ to $x$. (See Section~\ref{transp} for more details.)  
\end{definition}

We are now ready to define the notion of a reflector pair.
\begin{definition}
A pair $(\r,z)\in\Cp\times C(\Tb)$ is called a {\em{reflector pair}} if
$\rt,\zt>0$ and
\begin{align}
\rt(m)& = \sup_{x\in\Tb}\(\frac{1}{\zt(x)}K(m,x)\) \text{ for } m\in\Db, \label{def_refl_rt}\\
\zt(x)& = \sup_{m\in\Db}\(\frac{1}{\rt(m)}K(m,x)\) \text{ for } x\in\Tb. \label{refl_zt}
\end{align}
\end{definition}
Here we used the definitions of $\rt$ and $\zt$ from Definitions~\ref{def_rt} and \ref{def_zt}.
Note that the suprema on the right hand sides are in fact attained.
Also note that for a reflector pair $(\r,z)$, we have
\begin{equation}\label{def_refl}
\rt(m)\, \zt(x)\geq K(m,x)  \quad\quad {\text{ for all }}m\in\Db, x\in\Tb.
\end{equation}

We will show later that reflector pairs are in fact concave (see Proposition~\ref{convex})
and that $\log\rt$ and $\log\zt$ are uniformly Lipschitz continuous (see  Proposition~\ref{lipschitz}). 

Finally, we define the reflector map, or 
ray-tracing map, associated to a reflector pair. Again,
the choice of terminology will become clear when we consider the problem from a geometric
viewpoint in Section~\ref{geometry_refl_map}.

\begin{definition} \label{def_refl_map}
Let $(\r,z)\in\Cp\times C(\Tb)$ be a reflector pair. Define its {\em{reflector map}},
or {\em{ray tracing map}}, 
as a set-valued map $\gamma\colon \Db\to\{\text{subsets of }\Tb\}$ via
\[
\gamma(m)=\{x\in\Tb\;\bigl|\; \rt(m)=\frac{1}{\zt(x)}K(m,x)\}\quad \text{ for }m\in\Db.  
\]
Clearly $\gamma(m)\neq\emptyset$ for any $m\in\Db$ by (\ref{def_refl_rt}). We will show later
that $\gamma(m)$ is in fact single-valued for almost all $m\in\Db$.  
(See Proposition~\ref{singlevalued}.) We may thus regard $\gamma$ as a transformation
$\gamma\colon\Db\to\Tb$.
\end{definition}

\section{Reflector pairs: Geometric viewpoint} \label{section_geometry}
In this section, we will investigate the definition of reflector pairs from a geometric
point of view. We show that the reflectors can be obtained as the boundary of certain
convex sets. These sets in turn are the intersections of a family of spheroids (for reflector 1)
and  paraboloids (for reflector 2), respectively.

In the first two sections, we will use the following notation, the pointwise analogue of Definitions~\ref{def_zt}
and \ref{def_rt}: For $x\in\mathbb R^2$ and $z>0$, and $m\in S^2$ and $\r>0$, write
\begin{align}
\zt&=\delta - \frac{z}{\ell^2-|x|^2}  \label{zt} \\
\rt&=-\delta + \frac{1}{2\r\, (m_z+1)}.\label{rt}
\end{align}
\subsection{The spheroids $S_{x,z}$}
We define first a family of spheroids, indexed by points $(x,z)$.
\begin{definition}
Let $x\in\Tb$ and $z\in\mathbb R$ such that $\zt>0$. Define the set
\begin{equation}
S_{x,z}=\{\rho\cdot m\;\bigl|\; m\in S^2, \rt=\frac{1}{\zt}{K(m,x)}\}. \label{S}
\end{equation}
Here we used the notation from~(\ref{zt})
and (\ref{rt}). 
\end{definition}

\begin{figure}
\begin{center}
\end{center}
\caption{LEFT: The set $S_{x,z}$ is a spheroid with foci at the origin $O$ and at $(x,z)$.
$\bar{S}_{x,z}$ is the convex set bounded by $S_{x,z}$. The sketch shows a cross section
through a plane containing the axis defined by the foci. The sets $S_{x,z}$
and  $\bar{S}_{x,z}$ are obtained by rotating the gray ellipse around the axis through the foci as
indicated.
RIGHT: The set $P_{x,z}$ is a paraboloid with focus at $\r\cdot m$ and axis
parallel to the $z-$axis.
$\bar{P}_{\r\cdot m}$ is the convex set bounded by $P_{\r\cdot m}$. 
The sketch shows a cross section
through a plane containing the axis and the focus. The sets  $P_{\r\cdot m}$
and $\bar{P}_{\r\cdot m}$  are obtained by rotating the gray parabola around the axis as
indicated.
}\label{fig_spheroid}\label{fig_paraboloid}
\end{figure}

\begin{lemma}
Using the notation from the previous definition, the condition
$\zt>0$ implies  $\ell-z>|(x,z)|$.  (Here $|.|$ denotes the standard Euclidean vector norm.)
The set $S_{x,z}$ is given by the following equation for $\rho>0$, $m\in S^2$: 
\begin{equation}
     \rho+|(x,z)-\rho\cdot m|=\ell-z.\label{cond1}
\end{equation}
Geometrically, $S_{x,z}$ is a spheroid whose 
foci are the origin $O=(0,0,0)$
and the point $(x,z)$. (See Figure~\ref{fig_spheroid}.)
\end{lemma}
\begin{proof}
Note first that the condition $\zt>0$ implies
\[
z=(\delta-\zt)(\ell^2-|x|^2)<\d (\ell^2-|x|^2)=\frac{\ell}{2}-\frac{|x|^2}{2\ell},
\]
and so
\[
\ell-z>\frac{\ell}{2}+\frac{|x|^2}{2\ell}>0.
\]
Also,
\[
(\ell-z)^2=\ell^2+z^2-2\ell z>z^2+|x|^2=|(x,z)|^2.
\]
It follows that
\[
\ell-z>|(x,z)|.
\]
A straightforward but lengthy algebraic computation now yields that the condition in (\ref{S}) 
is equivalent to 
\begin{equation}
(\ell-z-\r)^2=|(x,z)-\r\cdot m|^2.\label{help1}
\end{equation}
Here $\r=(2(m_z+1)(\rt+\d))^{-1}>0$, as $\rt>0$.
We have
\begin{align*}
\ell-z>|(x,z)|\geq \r-|(x,z)-\r\cdot m|=\r-|\ell-z-\r|,
\end{align*}
and thus $\ell-z-\r>-|\ell-z-\r|$. It follows that
\[
 \ell-z-\r=|\ell-z-\r|.
\]
This together with (\ref{help1}) gives equation (\ref{cond1}). 
Now consider the geometric content of equation (\ref{cond1}).
We can immediately read off that the set $S_{x,z}$ consists of all points
that satisfy that the sum of the distances to the points $(0,0,0)$ and $(x,z)$
equals $\ell-z$. (See Figure~\ref{fig_spheroid}. Also note that 
$\ell-z>|(x,z)|$ implies that the set $S_{x,z}$ is nonempty.) This is by definition a spheroid.
 \end{proof}

\begin{definition}
Let as before $x\in\mathbb R^2$ and $z\in\mathbb R$ such that $\zt>0$. Denote by
$\bar{S}_{x,z}$ the closed convex set bounded by $S_{x,z}$. Thus $\bar{S}_{x,z}$
is given by
\begin{equation}
\bar{S}_{(x,z)}=\{\rho\cdot m\;\bigl|\; m\in S^2, \rt\geq\frac{1}{\zt}{K(m,x)}\}. \label{Sb}
\end{equation}
\end{definition}

\subsection{The paraboloids $P_{\r\cdot m}$}
Similarly to the previous section, we now define a family of paraboloids $P_{\r\cdot m}$.
\begin{definition}
Let $m\in \Db\subseteq S^2$ and $\rho>0$ such that $\rt>0$.
Define the set
\begin{equation}
P_{\r\cdot m}=\{(x,z)\;\bigl|\; x\in \mathbb R^2, \zt=\frac{1}{\rt}{K(m,x)}\}. \label{P}
\end{equation}
Here we used again the notation (\ref{zt}), (\ref{rt}).
\end{definition}

\begin{lemma}
With the notation of the previous definition, 
consider new  coordinates $(p,q)$ obtained by a shift by $\r\cdot m$, that is
 $(p,q)=(x,z)-\r\cdot m$. Then the set $P_{\r\cdot m}$ is given 
by the equation
\begin{equation}
    4 a q= |p|^2-4a^2 \label{cond2}
\end{equation}
where 
\[
a=\frac{2\r\d (1+mz)-1}{4\d}<0.
\]
Geometrically, $P_{\r\cdot m}$ is a paraboloid with focus at $\r\cdot m$ and axis parallel to the
$z-$axis. It opens in the directions of the negative $z-$axis. (See Figure~\ref{fig_paraboloid}.)
\end{lemma}
\begin{proof}
A straightforward computation yields that (\ref{P}) is equivalent to (\ref{cond2}).
Note that a paraboloid with focus at $(0,0,0)$ and focal parameter $2\alpha$ has the equation 
$4 \alpha z=|x|^2-4 \alpha^2$. Also note that $a<0$ follows from $\rt>0$. 
\end{proof}

\begin{definition}
Let as before $m\in \Db\subseteq S^2$ and $\rho>0$ such that $\rt>0$. 
Denote by
$\bar{P}_{\r\cdot m}$ the closed convex set bounded by $P_{\r\cdot m}$. Thus $\bar{P}_{\r\cdot m}$
is given by
\begin{equation}
\bar{P}_{\r\cdot m}=\{(x,z)\;\bigl|\; x\in\mathbb R^2, \zt\geq\frac{1}{\rt}{K(m,x)}\}. \label{Pb}
\end{equation}
\end{definition}

\subsection{Geometry of reflector pairs}
We now investigate the geometric content of the definition of reflector pairs.
For this, suppose that $(\r,z)\in \Cp\times C(\Tb)$ is a reflector pair.
Consider the two sets
\begin{align}
\bar{S}&=\bigcap_{x\in\Tb}\bar{S}_{x,z(x)},\label{Sint} \\
\bar{P}&=\bigcap_{m\in\Db}\bar{P}_{\r(m)m}\label{Pint}.
\end{align}
Note that none of the sets $\bar{S}_{x,z(x)}$ and $\bar{P}_{\r(m)m}$  are empty as in 
fact $\rt(m)>0$ for all $m\in\Db$
by definition of reflector pairs and $\zt(x)>0$ for all $x\in\Tb$. A compactness argument yields
that $\bar{P}$ and $\bar{S}$  are nonempty as well, and they are convex sets, 
since they are the intersections of convex sets. Moreover, the algebraic
representations (\ref{Sb}) and (\ref{Pb}) along with the definition of
reflector pairs immediately give rise to the following geometric facts:
\begin{proposition}
Let $(\r,z)\in \Cp\times C(\Tb)$ be a reflector pair.
\begin{enumerate}
\item The graph of $\r$,
\[
\Gamma_\r=\{\r(m)\cdot m\;\bigl|\; m\in \Db\},\\
\]
is contained in the boundary $\partial \bar{S}$ of the convex set $\bar{S}$ given in (\ref{Sint}).
In fact,
\[
 \Gamma_\r=\partial\bar{S}\cap(\mathbb R^+\cdot \Db), 
\]
where $\mathbb R^+\cdot \Db=\{r\cdot m\;\bigl|\;r\geq 0,  m\in \Db\}$ is a cone
with vertex at the origin and cross section $\Db$.
(See Figure~\ref{Fig_gammarho}.)
\item
The graph of $\r$,
\[
\Gamma_z=\{(x,z(x))\;\bigl|\; x\in \Tb\}
\]
is contained in the boundary $\partial \bar{P}$ of the convex set $\bar{P}$ given in (\ref{Pint}).
In fact,
\[
 \Gamma_z=\partial\bar{P}\cap(\Tb\times\mathbb R), 
\]
where $\Tb\times\mathbb R=\{(x,h)\;\bigl|\;x\in\Tb, h\in\mathbb R\}$ is the cylinder
with base $\Tb$. (See Figure~\ref{Fig_gammaz}.)
\end{enumerate}
\end{proposition}
We may thus think of the two reflectors as the ``envelopes'' of certain families of
spheroids and paraboloids, respectively.
An immediate corollary is also the following result:
\begin{proposition}\label{convex}
If $(\r,z)\in\Cp\times C(\Tb)$ is a reflector pair, then $z$ is a concave function
and $\r$ is the radial function of a convex set. In particular, both $\r$ and $z$ are 
locally Lipschitz continuous and 
almost everywhere differentiable. (Here ``almost everywhere'' refers to the standard
Lebesgue measures on $S^2$ and $\mathbb R^2$, respectively.)
\end{proposition}
\begin{proof}
The first part of the above statement follows immediately from the previous discussion.
The Lipschitz and differentiablility properties are standard results from
convexity theories, see \cite{schneider}, Theorem 1.5.1 and Section 1.7.
\end{proof}

\begin{figure}
\begin{center}
\end{center}
\caption{Geometry of reflector pairs. LEFT: The graph $\Gamma_\r$ of $\r(m)$ 
is obtained by intersecting the
boundary of the convex set $\bar{S}$ with the cone $\mathbb R^+\Db$.
RIGHT: The graph $\Gamma_z$ of $z(x)$ is obtained by intersecting the
boundary of the convex set $\bar{P}$ with the cylinder $\Tb\times\mathbb R$.
}\label{Fig_gammaz}\label{Fig_gammarho}
\end{figure}

\subsection{Geometry of the reflector map}\label{geometry_refl_map}
We now investigate the geometry of the reflector map.
For this, the following terminology is useful:
Suppose $(\r,z)\in\Cp\times C(\Tb)$ is a reflector pair. For $m\in\Db$,
if $\r(m)\cdot m\in \Gamma_\r\cap S_{x,z(x)}$, say that the spheroid
$ S_{x,z(x)}$ is supporting to the graph $\Gamma_\r$ at the point $\r(m)\cdot m$.

Similarly, for $x\in\Tb$, if $(x,z(x))\in\Gamma_z\cap P_{\r(m)\cdot m},$
say the paraboloid $P_{\r(m)\cdot m}$ is supporting to $\Gamma_z$ at
$(x,z(x))$.

Note that $ S_{x,z(x)}$ is supporting to $\Gamma_\r$ at the point $\r(m)\cdot m$
if and only if $P_{\r(m)\cdot m}$ is supporting to $\Gamma_z$ at
$(x,z(x))$. This is because both statements are equivalent to
\[
   \zt(x)\,\rt(m)=K(m,x)
\]
by the definitions (\ref{P}) and (\ref{S}).

We also have the following geometric lemma:
\begin{lemma}
Let $x_1,x_2\in\Tb$ be two distinct point: $x_1\neq x_2$.
Suppose the spheroids $S_{x_1,z(x_1)}$ and $S_{x_2,z(x_2)}$
are supporting to $\Gamma_\r$ at the same point $\r(m)\cdot m$.
Then $S_{x_1,z(x_1)}$ and $S_{x_2,z(x_2)}$ intersect tranversally
at $\r(m)\cdot m$.
\end{lemma}
\begin{proof}
Assume the contrary, that is, that $S_{x_1,z(x_1)}$ and $S_{x_2,z(x_2)}$ intersect
tangentially to each other. Since the two spheroids share the focus
$O$, it follows from basic properties of ellipsoids that the two
line segments $\overline{(x_1,z(x_1)), \r(m)\cdot m}$ and
$\overline{(x_2,z(x_2)), \r(m)\cdot m}$ are parallel. Thus the three 
points  $(x_1,z(x_1))$,$(x_2,z(x_2))$ and $\r(m)\cdot m$ are collinear.
On the axis through these three points, the point $\r(m)\cdot m$
cannot lie between the other two points; that is $\r(m)\cdot m$
must be one of the end points of the line segment defined by the three
points. 
But this contradicts that the paraboloid $P_{\r(m)\cdot m}$
contains both $(x_1,z(x_1))$ and $(x_2,z(x_2))$.
\end{proof}

Now consider the reflector map $\gamma$ associated with a reflector pair $(\r,z)$
as defined in Definition~\ref{def_refl_map}. In the language defined above,
we may now say that $\gamma(m)$ is the set of all points $x\in\Tb$
such that that the spheroid
$ S_{x,z(x)}$ is supporting to the graph $\Gamma_\r$ at the point $\r(m)\cdot m$.
Consider the case where for some $m\in\Db$, the set
$\gamma(m)$ contains more than one point, say $\{x_1,x_2\}\subseteq\gamma(m)$.
By the previous lemma, the two spheroids
$ S_{x_1,z(x_1)}$ and $ S_{x_2,z(x_2)}$ intersect transversally at 
$\rho(m)\cdot m$. Thus $\rho$ is not differentiable at $m$. But the set of such points
has measure zero by Proposition~\ref{convex}. We thus immediately have the following
result:
\begin{proposition}\label{singlevalued}
Let $(\r,z)\in\Cp\times C(\Tb)$ be a reflector pair. Then its reflector map $\gamma\colon\Db\to\Tb$
is almost everywhere single-valued. That is, the set of points $m$ where $\gamma(m)$ 
is {\em{not}} single valued has measure zero with respect to the standard measure on $S^2$
as a submanifold of the measure space $(\mathbb R^3,\mu)$, where $\mu$ is the standard
Lebesgue measure.
\end{proposition}

\begin{figure}
\begin{center}
\end{center}
\caption{Geometry of the reflector map: A ray emitted from the origin $O$
will be reflected to a ray traveling in the negative
$z-$direction labeled by $x=\gamma(x)$. See the text for details.
}\label{refl_sketch}
\end{figure}

\begin{remark}
We now justify the terminology of the reflector map from an optical point of
view. See Figure~\ref{refl_sketch} for the following considerations.
Let $(\r,z)\in\Cp\times C(\Tb)$ be a reflector pair. Let $m\in\Db$
such that $\gamma(m)$ is single valued. We show that under the 
geometric optics approximation, a ray emitted in the direction
$m\in\Db$ will be reflected off the first reflector $\Gamma_\r$
and then the second reflector $\Gamma_z$ in such a way that
the reflected ray is parallel to the negative $z-$axis and
that it intersects a plane perpendicular to the $z-$axis
in the point $x=\gamma(m)$.

Consider the reflection off reflector 1 first.
Since the spheroid $S_{x,z(x)}$ is tangential
to the reflector $\Gamma_\r$, the ray will be reflected
off $\Gamma_\r$ the same way it would be reflected off
 $S_{x,z(x)}$. By the geometrical properties of spheroids,
this means that the ray is reflected towards the focus $(x,z(x))$.
There, the ray will encounter $\Gamma_z$. It will be 
reflected the same way as it would be reflected off
the paraboloid $P_{\r(m)\cdot m}$, that is, in the 
direction of the negative $z-$axis.

Thus our definition of the reflector map is in agreement
with the physical law of reflection. In
the case when $\gamma(m)$ is multi-valued,
the first reflector has a singular point and a ray
will split up into a cone of light rays. These rays
will generate a set of directions whose projection
onto a plane perpendicular to the $z-$axis is $\gamma(m)$.
This is consistent with the physical phenomenon of
diffraction at singularities.
\end{remark}

The following statements about the reflector map
$\gamma(m)$ are analogous to Theorem 4.8 and
Lemma 4.9 in \cite{go04}. See this paper and
the the reference \cite{oliker_refl} for further details on the proofs.

\begin{theorem}
Let $\mathcal B$ denote the $\sigma$ algebra of
Borel sets on $\Tb$. Then for any subset $\tau\in\mathcal B$,
$\gamma^{-1}(\tau)$ is measurable relative to the standard
Lebesgue
measure of $\Db$. In addition, for any
non-negative locally integrable function $I$
on $\Db$, the function
\[
\mathcal L(\tau)=\int_{\gamma^{-1}(\tau)}I(m)\dm
\]
is a non-negative completely additive
measure on $\mathcal B$. (Here $\dm$ is the
standard measure on $S^2\subseteq \mathbb R^3$.)
 \end{theorem}

\begin{lemma}\label{lemma_change_var}
With the notation of the above Theorem, let $h$
be a continuous function on $\Tb$. Then we have
\[
\int_{\Tb}h(x)\mathcal L(dx)=\int_{\Db}h(\gamma(m))I(m)\dm.
\]
\end{lemma}

\section{The reflector problem and an equivalent constrained minimization problem}\label{section_equiv}
Let now $I$ be a nonnegative, integrable function on $\Db$, and $L$ be a nonnegative,
integrable function on $\Tb$, such that
\begin{equation}\label{pres_tot_energy}
\int_{\Db}I(m)\dm=\int_{\Tb}L(x)dx.
\end{equation}
We may interpret $I$ and $L$ as the intensity distribution functions of the light beams on the input
and output apertures, respectively. The above integral condition is simply (total)
energy conservation.

In this section of the paper, we now formulate the reflector
problem. More specifically, we can call this formulation a ``weak'' version
of the reflector problem since we do not require the input functions $I$ and $L$
to be differentiable, nor do we require the reflectors to be smooth surfaces.

We then formulate a second problem, which is an infinite dimensional linear programming
problem. One of the main results is that the two problems are in fact equivalent. This 
is stated and proved at the end of this section.

Let us first formulate the reflector problem:

\begin{problem} {\bf{(Reflector Problem)}}\label{probl1}
For given input and output intensities $I$ and $L$ satisfying (\ref{pres_tot_energy}),
find a pair $(z,\r)\in\Cp\times C(\Tb)$ that satisfies the following conditions:
\begin{enumerate}
\item $(z,\r)$ is a reflector pair
\item The reflector map $\gamma\colon\Db\to\Tb$ satisfies \label{loc_energy_pres}
\[
\int_{\gamma^{-1}(\tau)}I(m)\dm=\int_\tau L(x)dx
\]
for any Borel set $\tau\subseteq\Tb$.
\end{enumerate}
\end{problem}
This formulation builds on the geometrical 
interpretation of reflector maps as presented in Section~\ref{section_geometry}. 
Note that condition \ref{loc_energy_pres} is local energy conservation.

We have the following immediate corollary from Lemma~\ref{lemma_change_var}:
\begin{corollary}\label{changevar}
Let $(z,\r)$ be a solution to Problem~\ref{probl1}. Then we have
\[
    \int_{\Db} h(\gamma(m)) I(m)\dm= \int_\Tb h(x) L(x)dx
\]
for all functions $h\in C(\Tb)$.
\end{corollary}

Before we now formulate Problem~\ref{probl2}, we define the following
function space:
\begin{definition}
Define the set of {\em{admissible functions}} as
\begin{align*}
\Adm=\{(r,\z)\in C(\Db)\times C(\Tb)\;\bigl|\; r(m)+\z(x)&\geq \log K(m,x)\\
& \text{for all }m\in\Db, x\in\Tb\}.
\end{align*}
\end{definition}

\begin{problem}\label{probl2}
Minimize the functional
\[
\F(r,\z)=\int_\Db r(m) I(m)\dm+\int_\Tb \z(x)L(x)dx
\]
on the space $\Adm$.
\end{problem}
The two problems are equivalent, as expressed in the following theorem.

\begin{theorem} \label{thm_equiv}
Let $(\r,z)\in \Cp\times C(\Tb)$ be a reflector pair. Then $(\log\rt,\log\zt)\in\Adm$.
The following statements are equivalent:
\begin{enumerate}
\item\label{solp1} $(\r,z)$ solves the Reflector Problem~\ref{probl1}.
\item\label{solp2} $(\log\rt,\log\zt)$ minimizes the functional $\F$ on $\Adm$.
\end{enumerate}
\end{theorem}
\begin{proof}
The statement $(\log\rt,\log\zt)\in\Adm$ follows immediately from (\ref{def_refl}).
The proof of the equivalence of \ref{solp1} and \ref{solp2}
is analogous to those of Theorem 5.2 in \cite{go04}
and Theorem 3.4 in \cite{go03}. See also Theorem 1 in \cite{Gangbo/Mccann:95}.
Since this theorem is central to this paper, we give an outline of the
proof, omitting some of the technicalities, which can be filled in
with the above references.
\bigskip

$\ref{solp1}\Rightarrow\ref{solp2}:$ Suppose  $(\r,z)$ solves Problem~\ref{probl1}
with corresponding reflector map $\gamma$.
Let $(r,\z)\in\Adm$. Then for any 
$m\in\Db$ such that $\gamma(m)$  is single-valued, we have
\[
 r(m)+\z(\gamma(m))\geq\log K(m,\gamma(m))=\log\rt(m)+\log\zt(\gamma(m)).
\]
This yields
\begin{align*}
   \int_\Db r(m)I(m)\dm+&\int_\Db\z(\gamma(m))I(m)\dm \\
& \geq\int_\Db \log\rt(m)I(m)\dm+\int_\Db\log\zt(\gamma(m))I(m)\dm.
\end{align*}
Now using Corollary~\ref{changevar} gives
$
  \F(r,\z)\geq\F(\log\zt,\log\rt).
$
\bigskip

$\ref{solp2}\Rightarrow\ref{solp1}:$ The main idea is that the Euler-Lagrange
equations of minimizing $\F$ are equivalent to the equality
\[
    \int_{\Db} h(\gamma(m)) I(m)\dm= \int_\Tb h(x) L(x)dx
\]
for all functions $h\in C(\Tb)$. This implies that $(\r,z)$ solves Problem~\ref{probl1}.

Let thus $h\in C(\Tb)$. Let $\e>0$ be a small parameter.
To bring out
the main ideas, we present a formal calculation, assuming expansions in $\e$ are valid.
A completely rigorous treatment is 
possible; indeed the proof in \cite{go04} can easily be modified to the
problem at hand. 

Define perturbations $(r_\e, \z_\e)\in\Adm$
of $(\log\rt,\log\zt)$ via
\begin{align*}
\z_\e(x)&=\log\zt(x)+\e\, h(x) \quad{\text{for }}x\in\Tb\\ 
r_\e(m)&=\sup_{x\in\Tb}\(-\z_\e(x)+\log K(m,x)\) \quad{\text{for }}m\in\Db.
\end{align*}
Let now $x_\e$ be a point where the supremum in the definition of
$r_\e(m)$ is attained. Expanding $x_\e$ in $\e$ yields
\[
  x_\e=\gamma(m)+\O(\e).
\]
Thus again an expansion in $\e$ gives
\begin{align*}
r_\e(m)&=-\z_\e(x_\e)+\log K(m,x_\e)\\
       &=-\log\zt(\gamma(m))+\log K(m,\gamma(m))-\e\, h(\gamma(m))+\O(\e^2)\\
       &=\log\rt(m)-\e\, h(\gamma(m))+\O(\e^2).
\end{align*}
Thus, using the fact that $(\log\rt,\log\zt)$ minimizes $\F$,
\begin{align*}
0&=\frac{d}{d\e}\Bigl|_{\e=0}\F(r_\e,\z_\e)
 =-\int_{\Db} h(\gamma(m)) I(m)\dm+ \int_\Tb h(x) L(x)dx.
\end{align*}
This completes the sketch of the proof.
\end{proof}

\section{Existence of solutions}\label{section_existence}
Note that Theorem~\ref{thm_equiv}, while showing that Problems~\ref{probl1}
and \ref{probl2} are equivalent, does not state that solutions exist.
We prove this in the present section.

In the following, fix some point $m^*\in\Db$.
\begin{definition}
Set
\begin{align*}
\Admr=\{(\log\rt,\log\zt)&\in C(\Db)\times C(\Tb)\;\bigl|\; \\
& (\r,z) {\text{ is a reflector pair with }}\log\rt(m^*)=0\}.  
\end{align*}
Note that $\Admr\subseteq\Adm$.
\end{definition} 
\begin{proposition} \label{lipschitz}
The family of pairs of functions $\Admr$ is uniformly Lip\-schitz continuous
in each entry.
That is, there are constants $K_1,K_2>0$ such that
\begin{align*}
&|\z(x_1)-\z(x_2)|\leq K_1\cdot |x_1-x_2|\\
&|r(m_1)-r(m_2)|\leq K_2\,\cdot d(m_1,m_2)
\end{align*}
for all $(r,\z)\in\Admr$, $x_1,x_2\in\Tb$, $m_1,m_2\in\Db$.
Here $d(m_1,m_2)={\rm{dist}}_{S^2}(m_1,m_2)$ is the intrinsic distance on $S^2$.
The constants $K_1,K_2$ only depend on $\Db,\Tb$; explicitly, we have
\begin{align*}
K_1&=\max_{m\in\Db,x\in\Tb}|\nabla_x\log K(m,x)|,\\
K_2&=\max_{m\in\Db,x\in\Tb}|\nabla_m\log K(m,x)|.
\end{align*}
(Here $\nabla_m$ denotes the gradient with respect to the variable $m$
on the sphere $S^2$, and $\nabla_x$ is the gradient with
respect to $x$.)
\end{proposition}
\begin{proof}
We prove the inequality for $\z$; the proof of the other inequality
is completely analogous. Let $x_1,x_2\in\Tb$ and $(r,\z)\in\Admr$.
Assume $\z(x_2)\leq\z(x_1)$, otherwise relabel $x_1,x_2$.
Let $m_1\in\Db$ be such that $\z(x_1)+r(m_1)=\log K(m_1,x_1)$.
(Such an $m_1$ exists by (\ref{refl_zt}).) 
Then
\begin{align*}
0\leq \z(x_1)-\z(x_2)&=-\z(x_2)-r(m_1)+\log K(m_1,x_1)\\
 &\leq \log K(m_1,x_1)-\log K(m_1,x_2)\\
& \leq \(\max_{m\in\Db,x\in\Tb}|\nabla_x\log K(m,x)|\)\cdot |x_1-x_2|.
\end{align*}
\end{proof}

\begin{proposition}\label{bdd}
Functions in $\Admr$ are uniformly bounded. Specifically, we have
\begin{align*}
&|r(m)|\leq K_1\,\diam(\Db)\\
&|\z(x)|\leq \max_{x'\in\Tb}|\log K(m^*,x')|+K_2\,\diam(\Tb). 
\end{align*}
for all $(r,\z)\in\Admr$, $x\in\Tb$, $m\in\Db$. Here $K_1,K_2$ are as in
Proposition~\ref{lipschitz}, and $\diam(\Tb)=\max_{x_1,x_2\in\Tb}|x_1-x_2|$
and $\diam(\Db)=\max_{m_1,m_2\in\Db}d(m_1,m_2)$ are the diameters of $\Tb$ and $\Db$,
respectively.
\end{proposition}
\begin{proof}
Let $(r,\z)\in\Admr$, $x\in\Tb$, $m\in\Db$. By Proposition~\ref{lipschitz}, we have
\[
|r(m)|\leq |r(m^*)|+K_1\, d(m,m^*)\leq K_1\,\diam(\Db).
\]
Also, for some $x^*\in\gamma(m^*)$, we have $\z(x^*)=\log K(m^*,x^*)$, and thus again
by Proposition~\ref{lipschitz},
\begin{align*}
|\z(x)|\leq |\z(x^*)|+K_2\, |x-x^*|
  \leq \max_{x'\in\Tb}|\log K(m^*,x')|+K_2\,\diam(\Tb).
\end{align*}
\end{proof}

We can now prove the existence of a solution to Problem~\ref{probl2}.
With Theorem~\ref{thm_equiv}, this immediately implies the
existence of a solution to the reflector problem.

\begin{theorem}
The functional $\F$ attains its minimum on $\Adm$. Moreover, this minimum
is actually attained at some $(r,\z)\in\Admr$.
\end{theorem}
\begin{corollary}\label{existence}
The Reflector Problem~\ref{probl1} has a solution.
\end{corollary}
\begin{proof}
By Propositions~\ref{lipschitz} and \ref{bdd} and the Arzel\`a-Ascoli Theorem,
the functional $\F$ attains its minimum on $\Admr$.

We now show that this minimum is also the minimum of $\F$ on the larger set $\Adm$.
For this, let $(r,\z)\in\Adm$. Define
\begin{align*}
r^*(m)=\sup_{x\in\Tb}\{-\z(x)+\log K(m,x) \}\quad \text{ for }m\in\Db\\
\z^*(x)=\sup_{m\in\Db}\{-r^*(m)+\log K(m,x) \}\quad \text{ for }x\in\Tb
\end{align*}
It follows that
\[
 r^*(m)\leq r(m),\quad\quad\quad \z^*(x)\leq\z(x)
\]
for all $m\in\Db$, $x\in\Tb$. Note that $\z^*(x)\leq\z(x)$
implies
\[
  \sup_{x\in\Tb}\{-\z^*(x)+\log K(m,x) \}\geq \sup_{x\in\Tb}\{-\z(x)+\log K(m,x) \}=r^*(m).
\]
But $-\z^*(x)+\log K(m,x)\leq r^*(m)$ for all $x\in\Tb,m\in\Db$ implies
\[
  \sup_{x\in\Tb}\{-\z^*(x)+\log K(m,x) \}\leq r^*(m).
\]
It follows that
\[
r^*(m)=\sup_{x\in\Tb}\{-\z^*(x)+\log K(m,x) \}\quad \text{ for }m\in\Db.
\]
Thus $(r^*,\z^*)=(\log\rt,\log\zt)$ for some reflector pair
$(\r,z)$. 
Note that
\[
(r^*-r^*(m^*), \z^*+r^*(m^*))\in\Admr,
\]
where the left hand pair denotes functions shifted by the constants
$\pm r^*(m^*)$. 
Using the above, we now have
\begin{align*}
\F(r^*-r^*(m^*), \z^*+r^*(m^*))&=\F(r^*,\z^*)=\int_\Db r^*(m)I(m)\dm+\int_\Tb \z^*(x)L(x)dx\\
&\leq \int_\Db r(m)I(m)\dm+\int_\Tb \z(x)L(x)dx= \F(r,\z).
\end{align*}
This shows that indeed the minimum of $\F$ on $\Adm$ is attained on $\Admr$.
The corollary is an immediate consequence of Theorem~\ref{thm_equiv}.
\end{proof}

\section{A uniqueness result and connection to an optimal transportation problem}\label{transp}
In this section, we show the connection to an optimal transportation problem.
This is again quite analogous to the reflector design problems in \cite{go04}
and \cite{go03}. This connection allows us to formulate a uniqueness result
for the Reflector Problem~\ref{probl1}.

For the formulation of the problem, we need the concept of a {\em{plan}}
in this context:
\begin{definition}
A {\em{plan}} is a map $P\colon\Db\to\Tb$ that is measure preserving,
that is, we have
\[
    \int_{\Db} h(P(m)) I(m)\dm= \int_\Tb h(x) L(x)dx
\]
for any function $h\in C(\Tb)$.
\end{definition}
A plan $P$ needs only be defined almost everywhere on $\Db$.

The optimal transportation problem associated with intensities
$I(m)$ and $L(x)$ satisfying energy conservation (\ref{pres_tot_energy})
and with cost function $\log K(m,x)$ is the following:
\begin{problem}
Maximize the transportation cost
\begin{equation}\label{transp_cost}
    P\mapsto C(P)=\int_{\Db}\log K(m,P(m))I(m)\dm   
\end{equation}
among all plans $P.$
\end{problem}

This problem is again solved by the reflector map
of a solution to the Reflector Problem~\ref{probl1}.
In fact, we have the following theorem:
\begin{theorem}\label{thm_transp_cost}
Let $(z,\r)$ be a solution to the Reflector Problem~\ref{probl1}.
Then the corresponding reflector map $\gamma$ is a plan,
and it maximizes the transportation cost (\ref{transp_cost})
among all plans. Any other cost maximizing plan is equal
to $\gamma$ almost everywhere on 
$\supp I\setminus\{m\in\Db\,\bigl|\, I(m)=0\}$.
\end{theorem}
\begin{proof}
Let $(z,\r)$ be a solution to the Reflector Problem~\ref{probl1}
with corresponding reflector map $\gamma$. 
By Corollary~\ref{changevar}, $\gamma$ is in fact a plan.
Now let $P$ be another plan. Then
\[
\log\rt(m)+\log\zt(P(m))\geq \log K(m,P(m))
\]
for almost all $m\in\Db$, and equality holds iff $P(m)=\gamma(m)$.
Thus
\begin{align*}
C(P)&=\int_{\Db}\log K(m,P(m))I(m)\dm\\
  &\leq \int_\Db\log\rt(m)I(m)\dm+\int_\Db\log\zt(P(m))I(m)\dm\\
  &= \int_\Db\log\rt(m)I(m)\dm+\int_\Tb\log\zt(x)L(x)dx\\
  &=\int_\Db\log\rt(m)I(m)\dm+\int_\Db\log\zt(\gamma(m))I(m)\dm\\
  &=\int_{\Db}\log K(m,\gamma(m))=C(\gamma).
\end{align*}
Thus $\gamma$ indeed maximizes the transportation cost among all plans. 
Moreover, if equality holds in the above estimate, then
$\gamma(m)=P(m)$ or $I(m)=0$ for almost all $m\in\Db$.
\end{proof}

We can now state the following theorem, summarizing the 
main result on the existence of solutions
to the Reflector Problem~\ref{probl1} with an additional uniqueness result:

\begin{theorem}\label{thm_ex_unique} (Existence and uniqueness for solutions to the Reflector Problem~\ref{probl1})
There exist solutions to the Reflector Problem~\ref{probl1}. If $(\r,z)$ is a solution,
then both $\r$ and $z$ are Lipschitz continuous.
The corresponding reflector map $\gamma\colon\Db\to\Tb$
is single valued almost everywhere on $\Db$. If $(\r,z)$
and $(\r',z')$ are two solutions with reflector maps $\gamma$ 
and $\gamma'$, respectively, then $\gamma(m)=\gamma'(m)$
for almost all $m\in\supp I\setminus\{m\in\Db\,\bigl|\, I(m)=0\}$.
\end{theorem}
\begin{proof}
The existence of solutions to Problem~\ref{probl1} 
was already obtained in Corollary~\ref{existence}. If now
$\gamma$ and $\gamma'$ are two reflector maps corresponding 
to two solutions, then both maximize the transportation cost $C$
by Theorem~\ref{thm_transp_cost}, and hence we have
$\gamma(x)=\gamma(x')$
for almost all $\supp I\setminus\{m\in\Db\,\bigl|\, I(m)=0\}$.
\end{proof}
There is a number of open questions for further investigations. 
For instance, we have assumed certain constraints
on the aperture $\Db$, in particular 
$(0,0,-1)\notin\Db$. From the physical intuition
about the problem, these constraints appear to be unnecessary. It would
be interesting to extend the theory to these cases as well.
Furthermore, a further exploration of the 
regularity of 
solutions $(\r,z)$ depending on the intensities
$I$ and $L$ would be desirable. It is expected that some
of the currently rapidly growing research on optimal
transportation may carry over here. (See 
the recent survey \cite{villani} and the extremely 
extensive bibliography cited there.)
 
\bigskip
{\bf{Acknowledgement:}} The author would like to thank V. Oliker
for useful discussions.

\nocite{wang96}
\bibliographystyle{plain}
\bibliography{bibfile}

\end{document}